\documentclass{amsart}[10pt]

\usepackage{epsfig}
\usepackage{graphics}
\usepackage{amsfonts}
\usepackage{amsmath}
\usepackage{latexsym}

\newtheorem{theorem}{Theorem}
\newtheorem{lemma}{Lemma}
\newtheorem{proposition}{Proposition}

\newtheorem{remark}{Remark} 

\begin{document}

\newcommand\ho{H\"older }
\newcommand\ep{\varepsilon}

\title[Inside singularity sets of measures]{Inside singularity sets of
random Gibbs measures}

%\runningtitle{ee}

\author{Julien Barral and St\'ephane Seuret}

\address[]{INRIA Rocquencourt, \'Equipe ``Complex'', Domaine de
Voluceau Rocquencourt, 78153 Le Chesnay cedex, France}

%%%%%%%%%%%%%%%%%%%%%%%%%%%%%%%
\begin{abstract}
We evaluate the scale at which the multifractal structure of some
random Gibbs measures becomes discernible. The value of this scale is
obtained through what we call the growth speed in H\"older singularity
sets of a Borel measure. This growth speed yields new information on
the multifractal behavior of the rescaled copies involved in the
structure of statistically self-similar Gibbs measures. Our results
are useful to understand the multifractal nature of various
heterogeneous jump processes.

%Results are given for random Gibbs measures.
\end{abstract}
%%%%%%%%%%%%%%%%%%%%%%%%%%%%%%%
%%%%%%%%%%%%%%%%%%%%%%%%%%%%%
\keywords{Random Gibbs measures; Self-similarity; Large Deviations;
Hausdorff dimension; Fractals}

\maketitle

%%%%%%%%%%%%%%%%%%%%%%%%%%%%%%%%%
%\ams{60G57|60F10|28A78}{28A80} 
% insert the primary Maths Subject Classification number in the first
% bracket
% and the secondary ams number(s) in the second bracket
% e.g. \ams{60E20}{49G03;49F10}

%%%%%%%%%%%%%%%%%%%%%%%%%%%%%%%%%%%%%%%%%%%%%%%%

%%%%%%%%%%%%%%%%%%%%%%%%%%%%%%%%%%%%%%%%%%%%%%%%
%%%%%%%%%%%%%%%%%%%%%%%%%%%%%%%%%%%%%%%%%%%%%%%%
%%%%%%%%%%%%%%%%%%%%%%%%%%%%%%%%%%%%%%%%%%%%%%%%
%%%%%%%%%%%%%%%%%%%%%%%%%%%%%%%%%%%%%%%%%%%%%%%%
\section{Introduction}
\label{intro}

Contrary to what happens with monofractal measures (for instance
uniform measures on regular Cantor sets), multifractal measures
exhibit simultaneously several different behaviors at small scales. It
is natural to question from which scale the multifractal structure of
these measures becomes discernible and remains stable. This paper
introduces a notion which provides a way to examinate the value of
this critical scale. This notion, that we call growth speed in
singularity sets, is naturally related with multifractal measures. In
the following, we define and study the growth speed in singularity
sets for a class of statistically self-similar measures which includes
random Gibbs measures. This work requires refinements of the known
theoretical results on the multifractal nature of these
measures. Finally, we obtain rigorous estimates of the error made when
approximating the asymptotic local behavior of the measure by
observing it at a fine but fixed grid.

\smallskip

Before making precise all these notions, let us explain what one of
our main motivations was.  The new multifractal properties we point
out in this paper are naturally involved in the small-scale structure
analysis of some jump processes recently considered in
\cite{MOIBARRALNU,NUGEN,LEVY}. Typical examples of such heterogeneous
jump processes are L\'evy processes in multifractal time. Performing a
multifractal time change in irregular processes is a natural idea when
trying to build multi-parameter processes
\cite{M2,BaMu,Riedi}. Indeed, such processes yield multifractal
objects with an interesting structure, that may be more realistic than
classical homogeneous jump processes (for instance like L\'evy
processes) for the purpose of modeling multifractal discontinuous
phenomena (Internet traffic \cite{JLVRIEDI}, variations of financial
prices \cite{M2}). Another relevant property of these processes is
that they provide new illustrations of multifractal formalisms
\cite{HALSEY,BRMICHPEY,OLSEN,MOIBARRALNU}.

 Our results provide tools to study these processes. Indeed, the
multifrac\-tal analysis of heterogeneous jump processes in
\cite{MOIBARRALNU,NUGEN,LEVY} requires to deepen our knowledge
regarding statistically self-similar singular measures generated by
multiplicative processes. The fact that these measures are locally
equivalent to a rescaled copy of themselves is exploited in a new
direction using the notion of growth speed in the H\"older singularity
sets of these copies. The growth speed yields new insights on the
structure of the process, which are more precise than those obtained
by only considering individually these copies as the same
probabilistic object. In particular, it provides a new quantitative
way of distinguishing two well-known families of statistically
self-similar singular measures, namely the random Gibbs measures
\cite{KAKI} and the independent random cascades, like Mandelbrot
canonical cascades \cite{M}.  This paper focuses on random Gibbs
measures, the case of the Mandelbrot canonical cascades is very
different and treated in~\cite{REFCAS}.

\smallskip
The multifractal structure of random Gibbs measures has been
extensively studied
(\cite{HALSEY,Rand,Kifer,FAN2,PESINWEISS,BCM,FeOl}). This topic is
concerned with the size estimation of the H\"older singularity sets of
such a measure $\mu$. These sets are defined as the level sets of the
pointwise H\"older exponent $\lim_{r\to 0^+}\frac{\log
\mu(B(t,r)}{\log (r)}$. The sizes of H\"older singularity sets are
measured through their Hausdorff (or packing) dimension. It can be
shown that these dimensions are obtained thanks to the Legendre
transform of a kind of free energy function $\tau_\mu$ related to
$\mu$. More precisely, let $b$ be an integer $\ge 2$ and $\mathcal{A}$
the alphabet $\{0,\dots,b-1\}$.  Suppose that we are working on the
symbolic space $\mathbb{A}=\mathcal{A}^{\mathbb{N}^*}$ endowed with
the product topology and the one-sided shift transformation
$\sigma$. If $w\in\mathcal{A}^*=\bigcup_{n\ge 1}\mathcal{A}^n$, the
$n$ step cylinder about $w$ in $\mathbb{A}$ is denoted by $[w]$. The
measures we are interested in are associated with some (random)
H\"older potentiel and the dynamical system $(\mathbb{A},\sigma)$.

The function $\tau_\mu$ considered in the multifractal formalism for
measures in \cite{HALSEY,BRMICHPEY} is obtained as follows: For every
$ q\in \mathbb{R}$, let
\begin{equation}
\label{tau}
\forall \,j\geq 1, \ \tau_{\mu,j}(q)=-\frac{1}{j}\log_b\,
\sum_{w\in\mathcal{A}^j}\mu ([w])^q \ \mbox{ and } \
\tau_{\mu}(q)=\liminf_{j\to\infty}\tau_{\mu,j}(q).
\end{equation}
The Legendre transform of $\tau_\mu$ at $\alpha>0$ is then
 $\tau_\mu^*(\alpha):=\inf_{q\in\mathbb{R}}\alpha q-\tau_\mu(q)$. Then
 the H\"older singularity set of level $\alpha>0$ is defined as
$$
E^\mu_\alpha=\Big \{t\in \mathbb{A}:\ \lim_{n\to\infty}\frac{\log_b
\mu\big ([t|n]\big )}{n}=\alpha\Big \}
$$ 
($t|n$ stands for $t_1\cdots t_n$). The Gibbs measures we consider
obey the multifractal formalism in the sense that $\dim\,
E_\alpha^\mu=\tau_\mu^*(\alpha)$ when $\tau_\mu^*(\alpha)>0$.

This property is classically implied by the existence of a probability
measure $\mu_\alpha$ of the same nature as $\mu$ and such that
$\mu_\alpha$ is concentrated on $E_\alpha^\mu\cap
E_{\tau_\mu^*(\alpha)}^{\mu_\alpha}$. This measure $\mu_\alpha$ is called
an analyzing measure of $\mu$ at $\alpha$.

\smallskip

The existence of the measure $\mu_\alpha$ has another important
consequence regarding the possibility of measuring how the mass of
$\mu$ is distributed at a given large enough scale. Indeed, a direct
consequence of the multifractal formalism (\cite{R}) and the existence
of $\mu_\alpha$ is that for any $\varepsilon>0$ and $\alpha>0$ such
that $\tau_\mu^*(\alpha)>0$, one has
\begin{equation}\label{ld}
\lim_{j\to\infty}\frac{\log_b\#\big \{w\in\mathcal{A}^j:
b^{-j(\alpha+\varepsilon)}\le \mu ([w])\le
b^{-j(\alpha-\varepsilon)}\big\}}{j}=\tau_\mu^*(\alpha).
\end{equation}

The result we establish in this paper brings precisions on these sizes
estimates. We consider a refined version of the sets $E_\alpha(\mu)$
by considering, for any sequence $\varepsilon_n$ going down to 0, the
sets
$$
\widetilde E^\mu_{\alpha,p}=\Big \{t\in \mathbb{A}:\forall \ n\ge p,\
b^{-n(\alpha+\varepsilon_n)}\le \mu\big ([t|n])\big )\le
b^{-n(\alpha-\varepsilon_n)}\Big\},$$
\begin{equation}\label{wtE}
\mbox{ and } \ \ \ \widetilde E^\mu_\alpha=\bigcup_{p\ge 1}E^\mu_{\alpha,p}.
\end{equation}
It is possible to choose $(\varepsilon_n)_{n\ge 1}$ so that with
probability one, for all the exponents $\alpha$ such that
$\tau_\mu^*(\alpha)>0$, one has $ \mu_\alpha(\widetilde
E^\mu_\alpha)=\Vert \mu_\alpha\Vert =1$.

Since the sets sequence $ \widetilde E^\mu_{\alpha,p}$ is
non-decreasing and $\mu_\alpha (\widetilde E^\mu_\alpha)=1$, the
growth speed $ GS(\mu,\alpha)$ in $\widetilde E^\mu_{\alpha,p}$ can be
defined as the smallest value of $p$ for which the
$\mu_\alpha$-measure of $\widetilde E^\mu_{\alpha,p}$ reaches a
certain positive fraction $f\in (0,1)$ of the mass of $\mu_\alpha$,
that is the number
$$
GS(\mu,\alpha)=\inf\Big \{p: \mu_\alpha (\widetilde
E^\mu_{\alpha,p})\ge f\, \Vert \mu_\alpha\Vert \Big \}.
$$

Now for $n\ge 1$ and $\alpha> 0$ let
\begin{equation}\label{Nn}
 \mathcal{N}_n(\mu,\alpha)=\#\big \{w\in\mathcal{A}^n:
b^{-n(\alpha+\varepsilon_n)}\le \mu ([w])\le
b^{-n(\alpha-\varepsilon_n)}\big\}.
\end{equation}
Heuristically, one has 
$$GS(\mu,\alpha)\approx \inf\{p:\forall n\ge p,\
b^{n(\tau^*(\alpha)-\varepsilon_n)}\le \mathcal{N}_n(\mu,\alpha)\le
b^{n(\tau^*(\alpha)+\varepsilon_n)}\},$$
i.e. $ GS(\mu,\alpha)$
controls by above the smallest rank $p$ from which considering the
evaluation of $\mathcal{N}_n(\mu,\alpha)$ at any scale $b^{-n}$
smaller than $b^{-p}$ yields a correct representation of the
asymptotic behavior of $\mathcal{N}_n(\mu,\alpha)$.

\smallskip

Our results concern estimates of the growth speed of singularities
sets of copies of $\mu$ involved in the self-similarity property of
$\mu$. To illustrate our purpose, let us describe the model of
statistically self-similar measures we shall work with in the
sequel. We shall consider a natural random counterpart to {\it
quasi-Bernoulli} measures introduced in \cite{MICHON2,BRMICHPEY} and
mainly illustrated by deterministic Gibbs measures on $\mathbb{A}$. We
are inspired in particular by self-similar Riesz products and their
random version constructed with random phases (see \cite{FAN2} and
examples of Section~\ref{Bernoulli}).

%%%%%%%%%%%%%%%%%%%%%%%%%%%%%%%%%%%%%%
%%%%%%%%%%%%%%%%%%%%%%%%%%%%%%%%%%%%%%
\subsection{Quasi-Bernoulli independent random measure} 
$\ $

In the sequel $\equiv$ means equality in distribution.

A random probability measure $\mu=\mu(\omega)$ on $\mathbb{A}$ is said
to be a {\it quasi-Bernoulli independent random measure} if there
exists a constant $C>0$ and two sequences of random measures
$(\mu_j)_{j\ge 1}$ and $(\mu^{(j)})_{j\ge 1}$ such that for every
$j\ge 1$,

\smallskip
\noindent
$\bullet$ \!
{\bf (P1)} \! $\forall \, (v,w)\in \mathcal{A}^j\times\mathcal{A}^*$, $
\frac{1}{C}\mu_j([v])\mu^{(j)}([w])\le \mu([vw])\le C
\mu_j([v])\mu^{(j)}([w])$,

\smallskip
\noindent
$\bullet$ \!\!
{\bf (P2)} \!\!  for every $r\in \{0,\dots, b-1\}$, $0<\mbox{\rm{ess}}\
\mbox{\rm{inf}}\ \mu([r])\le \mbox{\rm{ess}}\,\mbox{\rm{sup}} \
\mu([r])<\infty$,

\smallskip
\noindent
$\bullet$ \!\!  {\bf (P3)} \!\! $\big (\mu^{(j)}([w])\big
)_{w\in \mathcal{A}^*}\equiv \big (\mu([w])\big )_{w\in
\mathcal{A}^*}$.  $\mu$ is also denoted $\mu^{(0)}$,

\smallskip
\noindent
$\bullet$ \!\! 
{\bf (P4)} \!\! $\sigma (\mu_j([v]):\! v\in
\mathcal{A}^j)$ and $\sigma (\mu^{(j)}([w]):\! w\in \mathcal{A}^*)\!$ are
independent.

\smallskip

The measures $\mu^{(j)}$ are the copies of $\mu$ mentioned in the
 paragraphs above.

\smallskip

%%%%%%%%%%%%%%%%%%%%%%%%%%%%%%%%%%%%%%%%%%%%%%%%%%%%%%%%%%
%%%%%%%%%%%%%%%%%%%%%%%%%%%%%%%%%%%%%%%%%%%%%%%%%%%%%%%%%%

%%%%%%%%%%%%%%%%%%%%%%%%%%%%%%%%%%%%%%%
%%%%%%%%%%%%%%%%%%%%%%%%%%%%%%%%%%%%%%%
%%%%%%%%%%%%%%%%%%%%%%%%%%%%%%%%%%%%%%%
\subsection{Controlling the growth speed in H\"older singularity sets of the $(\mu^{(j)})$'s}

Let $\mu$ be quasi-Bernoulli independent measure. For each copy
$\mu^{(j)}$ of $\mu$, the corresponding family of analyzing measures
$\mu^{(j)}_\alpha$ will be defined as $\mu_\alpha$ is defined for
$\mu$. The result we focus on is the asymptotic behavior~of
\begin{equation}\label{GS}
GS(\mu^{(j)},\alpha)=\inf\Big \{N: \mu^{(j)}_\alpha \big (\widetilde
E^{\mu^{(j)}}_{\alpha,N}\big )\ge f\, \Vert \mu^{(j)}_\alpha\Vert
\Big \} \mbox{ as $j\to\infty$.}
\end{equation}
For sake of simplicity, we give in this introduction a shorter version
of our main result (Theorem~\ref{quasi-Bernoulli_random_ren}).
\medskip

%%%%%%%%%%%%%%%%%%%%%%%%%%%%%%%%%%
\noindent
{\bf Theorem A.} {\em Suppose that $\tau_\mu$ is $C^2$. With
probability one, for all $\alpha>0$ such that $\tau_\mu^*(\alpha)>0$
there exists $\beta>0$ such that if $j$ is large enough, $
GS(\mu^{(j)},\alpha)\le \exp\sqrt{\beta \log j}$.}
%%%%%%%%%%%%%%%%%%%%%%%%%%%%%%%%%%

\medskip

Let us introduce the quantity
$$
GS'(\mu^{(j)},\alpha)=\inf \big \{p: \forall\ n\ge p,
b^{n(\tau_\mu^*(\alpha)-\varepsilon_n)}\le
\mathcal{N}_n(\mu^{(j)},\alpha)\le
b^{n(\tau_\mu^*(\alpha)+\varepsilon_n)}\big \}.
$$
Theorem A also implies a control of $\mathcal{N}_n(\mu^{(j)},\alpha)$
(recall (\ref{Nn})).  A stronger version
(Theorem~\ref{quasi-Bernoulli_random_dev}) of the following result is
going to be proved.

\medskip

%%%%%%%%%%%%%%%%%%%%%%
%%%%%%%%%%%%%%%%%%%%%%
\noindent
{\bf Theorem B.} {\em Suppose that $\tau_\mu$ is $C^2$. The same
conclusion as in Theorem A holds if $GS(\mu^{(j)},\alpha)$ is replaced
by $ GS'(\mu^{(j)},\alpha)$.}

%%%%%%%%%%%%%%%%%%%%%%
%%%%%%%%%%%%%%%%%%%%%%

\medskip

As claimed above, Theorems A and B indeed yields new information on
the multifractal structure of random Gibbs measures.

Section~\ref{principes} contains new definitions and two propositions
that are used in Section~\ref{Bernoulli} and \ref{secproof2} to state
and prove stronger versions of Theorems A and
B. Section~\ref{secproof1} contains the proof of results concerning
the speed of convergence of $\tau_{\mu,j}$ to $\tau_\mu$.

We end this introduction by giving an application of Theorem A.

%%%%%%%%%%%%%%%%%%%%%%%%%%%%%%%%%%%%%%%%%%%%
%%%%%%%%%%%%%%%%%%%%%%%%%%%%%%%%%%%%%%%%%%%%
\subsection{An application: The Hausdorff dimension of new limsup sets}

Let $\mu$ be a quasi-Bernoulli independent random measure as defined
previously and consider $\nu$, its projection on $[0,1]$. Examples of
jump processes of \cite{MOIBARRALNU,LEVY}~are
$$\sum_{j\ge 0} \ \sum_{0\le k\le b^j-1} \!\! j^{-2} \, \nu
([kb^{-j},(k+1)b^{-j}]) \, \delta_{kb^{-j}}\ \mbox{ and } \ \big (X\circ
\nu ([0,t])\big )_{0\le t\le 1},$$
 where $X$ is a L\'evy process.
Basically, if $\{x_n\}$ denotes the countable set of jump points of
such a process and $(\lambda_n)_{n\ge 1}$ is a sequence decreasing to
0 such that $\limsup_{n\to\infty}B(x_n,\lambda_n)=[0,1]$, the
multifractal nature of these processes is closely related to the
computation of the Hausdorff dimension of the sets defined for every
$\alpha> 0$, $\xi>1$ by
$$
K(\alpha,\xi) \, = \, \bigcap_{N\ge 1} \ \ \bigcup_{n\ge 1: \,
\lambda_n^{\alpha+\varepsilon_n} \leq \nu
([x_n-\lambda_n,x_n+\lambda_n]) \leq \lambda_n^{\alpha-\varepsilon_n}}
[x_n-\lambda_n^\xi,x_n+\lambda_n^\xi]
$$
for some sequence $(\varepsilon_n)$ converging to 0. The set
$K(\alpha,\xi)$ contains the points that are infinitely often close
to a jump point $x_n$ at rate $\xi$ relatively to $\lambda_n$, upon
the condition that $\nu ([x_n-\lambda_n,x_n+\lambda_n])\sim
\lambda_n^{\alpha}$. This last condition implies that $\nu$ has
roughly a H\"older exponent $\alpha$ at scale $\lambda_n$ around
$x_n$. One of the main results of \cite{MOIBARRALNU,MOIBARRALUBIQUITY}
(see also \cite{NOTE}) is the computation of the Hausdorff dimension
of $ K(\alpha,\xi)$. Under a suitable assumption on $(\lambda_n)$,
it is proved in \cite{MOIBARRALNU,MOIBARRALUBIQUITY} that, with
probability one, for all $\alpha$ such that $\tau_\mu^*(\alpha)>0$ and
all $\xi\ge 1$,
\begin{equation}\label{DIM}
\dim\, K(\alpha,\xi)={\tau_\mu^*(\alpha)}/{\xi},
\end{equation}
where $\dim$ stands for the Hausdorff dimension. This achievement is a
non-trivial generalization of what is referred to as ``ubiquity''
properties of the resonant system $\{(x_n,\lambda_n)\}$. Ubiquity
plays a role for instance in the description of exceptional sets
arising in the problem of small denominators and the physical
phenomenon of resonance \cite{ARNOLD,DOD}. In the classical result,
$\nu$ is equal to the monofractal Lebesgue measure, so $\alpha =1$,
the condition $\lambda_n^{\alpha+\varepsilon_n} \leq \nu
([x_n-\lambda_n,x_n+\lambda_n]) \leq \lambda_n^{\alpha-\varepsilon_n}$
is trivial, and $\dim\, K(1,\xi)=1/\xi$ (see \cite{DOD} for
instance). 

The fact that, by Theorem A, the growth speed $ GS(\mu^{(j)},\alpha)$
behaves like $o(j)$ as $j\to\infty$ is a crucial issue in constructing
a Cantor set of Hausdorff dimension ${\tau_\mu^*(\alpha)}/{\xi}$ in
$K(\alpha,\xi)$.
%%%%%%%%%%%%%%%%%%%%%%%%%%%%%%%%%%%%%%%%%%%%%%%%
%%%%%%%%%%%%%%%%%%%%%%%%%%%%%%%%%%%%%%%%%%%%%%%%
%%%%%%%%%%%%%%%%%%%%%%%%%%%%%%%%%%%%%%%%%%%%%%%%
%%%%%%%%%%%%%%%%%%%%%%%%%%%%%%%%%%%%%%%%%%%%%%%%
\section{Definitions, Growth speed in singularity sets}
\label{principes}
%%%%%%%%%%%%%%%%%%%%%%%%%%%%%%%%%%%%%%%%%%%%%%%%
%%%%%%%%%%%%%%%%%%%%%%%%%%%%%%%%%%%%%%%%%%%%%%%%

In the sequel, $(\Omega,\mathcal{B},\mathbb{P})$ denotes the
probability space on which the random variables of this paper are
defined.

%%%%%%%%%%%%%%%%%%%%%%%%%%%%%%%%%%%%%%%%%%%%%%%%
%%%%%%%%%%%%%%%%%%%%%%%%%%%%%%%%%%%%%%%%%%%%%%%%
\subsection{Measure of singularity
sets: a neighboring boxes condition}\label{princ} $\ $

Let $\mu$ and $m$ be two probability measures with supports equal to
$\mathbb{A}$. 

\smallskip

With any $w\in \mathcal{A}^n$ can be associated the integer $i(w)\in
\{0,1,\ldots, b^n-1\}$ such that the $b$-adic subinterval of $[0,1]$
naturally encoded by $w$ is $[i(w)b^{-n},(i(w)+1)b^{-n}]$
(alternatively $i(w)b^{-n}=\sum_{k=1}^nw_kb^{-k}$). Then, if $(v,w)\in
\mathcal{A}^n$, $\delta(v,w)$ stands for $|i(v)-i(w)|$. This defines
an integer valued distance on $\mathcal{A}^n$. This distance yields a
notion of neighbors for cylinders of the same generation. This notion
coincides with the natural one on $b$-adic subintervals of the same
generation in $[0,1]$. 
%This notion is needed because when proving
%(\ref{DIM}) we consider the projection $\nu$ on $[0,1]$ of the Gibbs
%measure on $\mathbb{A}$, and the control the asymptotic behavior of
%$\nu$ on $b$-adic intervals and their neighbors is needed.

\smallskip

Let $\widetilde\varepsilon=(\varepsilon_n)_{n\ge 0}$ be a positive
sequence, $N\ge 1$, and $\beta\ge 0$.

We consider a slight refinement of the sets introduced in (\ref{wtE}):
For $p\ge 1$,
$$
 E^{\mu}_{\beta,p}(N,\widetilde\varepsilon) =\left \{\! t\in
\mathbb{A}\!:\!
\begin{cases} \forall \, n\ge p, \,\forall\ \gamma\in\{-1,1\},\\ \forall w\in
\mathcal{A}^n,\, \delta(w,t|n)\le N, \end{cases}
\hspace{-2mm} \, b^{\gamma
n(\beta-\gamma\varepsilon_n)}\mu([w])^{\gamma}\le 1\!\right \}\!\!\!
$$
\begin{equation}
\label{defebm}
\mbox{ and } \ \ \  E^{\mu}_\beta(N,\widetilde\varepsilon)=\bigcup_{p\ge
1}E^{\mu}_{\beta,p}(N,\widetilde\varepsilon).
\end{equation}
This set contains the points $t$ for which, at each scale $n$ large
enough, the $\mu$-measures of the $2N+1$ neighbors of $[t|n]$ (for the
distance $\delta$) belong to $[b^{- n(\beta+\varepsilon_n)},
b^{-n(\beta-\varepsilon_n)}]$. Controlling the mass of these neighbors
is necessary in the proof of (\ref{DIM}) when $\mu$ is a
quasi-Bernoulli independent random measure.

For $n\ge 1$ and $\varepsilon,\eta>0$, let us define the quantity
\begin{equation}
\label{defsneej}
S^{N,\varepsilon,\eta}_n(m,\mu,\beta)= \sum_{\gamma\in
\{-1,1\}}b^{n(\beta-\gamma\varepsilon)\gamma\eta} \hspace{-2mm}
\sum_{v,w\in \mathcal{A}^n: \ \delta(v,w)\le N} \hspace{-2mm}
m([v])\mu([w])^{\gamma\eta}.
\end{equation}

%%%%%%%%%%%%%%%%%%%%%%%%%%%%%%%%%%%%%
\begin{proposition}\label{principe1}
Let $(\eta_n)_{n\ge 1}$ be a positive sequence. \\ If $\sum_{n\ge
1}S^{N,\varepsilon_n,\eta_n}_n(m,\mu,\beta)<+\infty$, then
$E^{\mu}_\beta(N,\widetilde\varepsilon )$ is of full
$m$-measure.
\end{proposition} 
%%%%%%%%%%%%%%%%%%%%%%%%%%%%%%%%%%%%%%%%
\begin{remark}
The same kind of conditions was used in \cite{BBP} to obtain a
    comparison between the box \cite{BRMICHPEY} and centered
    \cite{OLSEN} multifractal formalisms.
\end{remark}
%%%%%%%%%%%%%%%%%%%%%%%%%%%%%%%%%%%%%%%%

%%%%%%%%%%%%%%%%%%%%%%%%%%%%%%%%%%%%%
\begin{proof}
For $\gamma\in\{-1,1\}$ and $n\ge 1$, let us define
\begin{equation}\label{ii}
E^{\mu}_\beta(N,\varepsilon_n,\gamma)=\left\{t\in\mathbb{A}:\begin{cases}
  \forall\ w\in\mathcal{A}^n,\\ \delta(w,t|n)\le N, \!\!\!\!
\end{cases} b^{\gamma
n(\beta-\gamma\varepsilon_n)}\mu([w])^{\gamma}\le 1\right \}. 
\end{equation}
For $t\in\mathbb{A}$, if there exists (a necessarily unique)
$w\in\mathcal{A}^n$ such that $i(w)-i(t|n)=k$, this word $w$ is
denoted $w_k(t)$. For $\gamma\in \{-1,1\}$, let $S_{n,\gamma}
=\sum_{-N\le k\le N}m_k$ with
$$m_k=m \big (\big \{t\in\mathbb{A}: b^{\gamma
n(\beta-\gamma\varepsilon_n)}\mu(w_k(t))^{\gamma}>1 \big \}\big
).
$$ 
One clearly has
\begin{equation}\label{BC}
m\big ((E^{\mu}_\beta (N,\varepsilon_n,-1))^c\bigcup (E^{\mu}_\beta
(N,\widetilde\varepsilon_n,1))^c\big )\le S_{n,-1}+S_{n,1},
\end{equation}
Fix $\eta_n>0$ and $-N\le k\le N$.  Let $Y(t)$ be the random variable
which equals $b^{\gamma n(\beta-\gamma\varepsilon_n)\eta_n}
\mu([w_k(t)])^{\gamma \eta_n}$ if $w_k(t)$ exists, and  0 otherwise.  The
Markov inequality applied to $Y(t)$ with respect to $m$ yields $m_k \le \int
Y(t) dm(t)$. Since $Y$ is constant over each cylinder $[v]$ of
generation $n$, we get
\begin{eqnarray*}
 m_k \le\sum_{v,w \in \mathcal{A}^n:
i(w)-i(v)=k} b^{n(\beta-\gamma\varepsilon_n)\gamma \eta_n}m([v])
\mu([w])^{\gamma \eta_n}.
\end{eqnarray*}
Summing over $|k|\le N$ yields $S_{n,-1}+S_{n,1}\le
S^{N,\varepsilon_n,\eta_n}_n(m,\mu,\beta)$.  The conclusion follows
from (\ref{BC}) and from the Borel-Cantelli Lemma.
\end{proof}
%%%%%%%%%%%%%%%%%%%%%%%%%%%%%%%%%%%%%

%%%%%%%%%%%%%%%%%%%%%%%%%%%%%%%%%%%%%
%%%%%%%%%%%%%%%%%%%%%%%%%%%%%%%%%%%%%
\subsection{Growth speed in families of singularity sets}\label{princ2}

Let $\Lambda$ be a set of indexes, and $\Omega^*$ a measurable subset
of $\Omega$ of probability 1.  Some notations and technical
assumptions are needed to state the result.

%%%%%%%%%%%%%%%%%%%%%%%%%%%%%%%%%%%%%
\smallskip

$\bullet$ For every $\omega\in\Omega^*$, we consider two sequences of
families of measures $\Big (\{\mu^{(j)}_{\lambda}\}_{\lambda\in
\Lambda}\Big)_{j\ge 0}$ and $\Big(\{m^{(j)}_{\lambda}\}_{\lambda\in
\Lambda}\Big)_{j\ge 0}$ such that for every $j\ge 0$, the elements of
the families $\{\mu^{(j)}_{\lambda}\}_{\lambda\in \Lambda}$ and
$\{m^{(j)}_{\lambda}\}_{\lambda\in \Lambda}$ are probability measures
on $\mathbb{A}$. For $\nu\in\{\mu,m\}$, $\{\nu^{(0)}
_{\lambda}\}_{\lambda\in \Lambda}$ is written
$\{\nu_{\lambda}\}_{\lambda\in \Lambda}$.

%%%%%%%%%%%%%%%%%%%%%%%%%%%%%%%%%%%%%
\smallskip

$\bullet$ We consider an integer $N\ge 1$, and a positive sequence
$\widetilde\varepsilon=(\varepsilon_{n})_{n\ge 1}$, as well as a
family of positive numbers
$(\beta_\lambda)_{\lambda\in\Lambda}$. Then, remembering (\ref{ii})
let us consider for every $j\ge 0$ and $p\ge 1$ the sets
\begin{equation}
\label{deflevel}
E^{\mu^{(j)}_\lambda}_{\beta_\lambda,p}(
N,\widetilde\varepsilon)=\bigcap_{n\ge p
}E_{\beta_\lambda}^{\mu^{(j)}_\lambda}(N, \varepsilon_{n},-1)\cap
E_{\beta_\lambda}^{\mu^{(j)}_\lambda}(N, \varepsilon_{n},1).
\end{equation}

$\bullet$ The sets $\{ E^{\mu^{(j)}_\lambda}_{\beta_\lambda,p}(
N,\widetilde\varepsilon)\}_p$ form a non-decreasing sequence. One then
defines the growth speed of $ E^{\mu^{(j)}_\lambda}_{\beta_\lambda,p}
( N,\widetilde\varepsilon)$ as the quantity
\begin{equation}
\label{defn}
GS(m^{(j)}_\lambda,\mu^{(j)}_\lambda,\beta_\lambda,
N,\widetilde\varepsilon)= \inf \Big\{
p\ge 1 : \ m_\lambda^{(j)}\big
(E^{\mu^{(j)}_\lambda}_{\beta_\lambda,p}( N,\widetilde\varepsilon)\big
) \ge {1}/{2} \Big\}.
\end{equation}
This number, maybe infinite, is a measurement of the number $p$ of
generations needed for $E^{\mu^{(j)}_\lambda}_{\beta_\lambda,p}(
N,\widetilde\varepsilon)$ to recover a certain given fraction (here
chosen equal to 1/2) of the probability measure $m_\lambda^{(j)}$. We
assume that $m_\lambda^{(j)}$ is concentrated on $\lim_{p\rightarrow
+\infty} E^{\mu^{(j)}_\lambda}_{\beta_\lambda,p}
(N,\widetilde\varepsilon)$, so that
$GS(m^{(j)}_\lambda,\mu^{(j)}_\lambda, \beta_\lambda,N,
\widetilde\varepsilon)<\infty$.

%%%%%%%%%%%%%%%%%%%%%%%%%%%%%%%%%%%%%
\smallskip

$\bullet$ We assume that for every positive sequence $\widetilde\eta=
(\eta_j)_{j\ge 0}$, there exist 

\noindent - a random vector $V (\widetilde\eta) \in
\mathbb{R}_+^{\mathbb{N}}$, a sequence $(V^{(j)})_{j\ge 0}$ of
copies of $V (\widetilde\eta)$,

\noindent - a sequence $(\psi_j(\widetilde\eta))_{j\ge 0}$ such that for
$\mathbb{P}$-almost every $\omega\in\Omega^*$,
\begin{equation}
\label{bound}
\forall \, j\ge 0, \ \  \forall n\ge \psi_{j} (\widetilde\eta), \ \ 
V^{(j)}_n \ge \sup_{\lambda\in\Lambda}S_n^{N,\varepsilon_{n},
\eta_n}(m^{(j)}_\lambda,\mu^{(j)}_\lambda,\beta_\lambda),
\end{equation}
where $S_n^{N,\varepsilon_{n},
\eta_n}(m^{(j)}_\lambda,\mu^{(j)}_\lambda,\beta_\lambda)$ is defined
in (\ref{defsneej}). This provides us with a uniform control over
$\lambda \in \Lambda$ of the families of measures
$(m^{(j)}_\lambda,\mu^{(j)}_\lambda)_{j\ge 0}$.

%%%%%%%%%%%%%%%%%%%%%%%%%%%%%%%%%%%%%
%%%%%%%%%%%%%%%%%%%%%%%%%%%%%%%%%%%%%
\begin{proposition}[{\bf Uniform growth speed in singularity sets}] 
\label{ren}
Let $\widetilde \eta=(\eta_j)_{j\ge0} $ be a the sequence of positive
numbers.\\ Let $(\mathcal{S}_j)_{j\ge 0}$ be a sequence of integers
such that $\mathcal{S}_j\ge \psi_j(\widetilde\eta)$. Assume that
\begin{equation}
\label{C1}
\sum_{j\ge 0}\sum_{n\ge
\mathcal{S}_j} \mathbb{E} \Big (V_n(\widetilde\eta)\Big )<\infty.
\end{equation}
With probability one, for every $j$ large enough, for every
$\lambda\in \Lambda$, one has
$GS(m^{(j)}_\lambda,\mu^{(j)}_\lambda,\beta_\lambda,
N,\widetilde\varepsilon)\le \mathcal{S}_j$.
\end{proposition}
%%%%%%%%%%%%%%%%%%%%%%%%%%%%%%%%%%%%%
%%%%%%%%%%%%%%%%%%%%%%%%%%%%%%%%%%%%%
%%%%%%%%%%%%%%%%%%%%%%%%%%%%%%%%%%%%%
%%%%%%%%%%%%%%%%%%%%%%%%%%%%%%%%%%%%%
\begin{proof}
Fix $j\ge 1$. As shown in Proposition \ref{principe1}, for every $n\ge
1$ and every $\lambda\in\Lambda$, one can write
\begin{eqnarray*}
m^{(j)}_\lambda\Big (\big(E^{\mu^{(j)}_\lambda}_{\beta_\lambda}(N,
\varepsilon_{n},-1)\big )^c\cup \big
(E^{\mu^{(j)}_\lambda}_{\beta_\lambda}(N,
\varepsilon_{n},1)\big )^c\Big )
 \le
S_n^{N,\varepsilon_{n},\eta_n}(
m^{(j)}_\lambda,\mu^{(j)}_\lambda,\beta_\lambda).
\end{eqnarray*}
Thus, using (\ref{bound}), one gets
\begin{equation}
\label{ineg1}
m^{(j)}_\lambda \Big ( \bigcup_{n\ge \mathcal{S}_j} \big
(E^{\mu^{(j)}_\lambda}_{\beta_\lambda}(N, \varepsilon_{n},-1)\big
)^c\cup \big (E^{\mu^{(j)}_\lambda}_{\beta_\lambda}(N,
\varepsilon_{n},1)\big )^c\Big )\le \sum_{n\ge
\mathcal{S}_j}V^{(j)}_n.\!\!
\end{equation}
Now (\ref{C1}) yields
\begin{eqnarray*}
\sum_{j\ge 1}\mathbb{P}\Big( \sum_{n\ge \mathcal{S}_j}V^{(j)}_n\ge
{1}/{2}\Big ) \le 2\sum_{j\ge 1}\mathbb{E}\Big (\sum_{n\ge
\mathcal{S}_j}V^{(j)}_n \Big )<\infty.
\end{eqnarray*}
Thus, with probability one, $ \sum_{n\ge \mathcal{S}_j}V^{(j)}_n<
{1}/{2} $ for every $j$ large enough. This, combined with
(\ref{deflevel}), (\ref{ineg1}) and (\ref{defn}), implies that, with
probability one, for all $j$ large enough, for every $\lambda\in
\Lambda$, $GS(m^{(j)}_\lambda,\mu^{(j)}_\lambda,\beta_\lambda,
N,\widetilde\varepsilon)\le \mathcal{S}_j$.
\end{proof}
%%%%%%%%%%%%%%%%%%%%%%%%%%%%%%%%%%%%%%

%%%%%%%%%%%%%%%%%%%%%%%%%%%%%%%%%%%%%%
%%%%%%%%%%%%%%%%%%%%%%%%%%%%%%%%%%%%%%
%%%%%%%%%%%%%%%%%%%%%%%%%%%%%%%%%%%%%%
%%%%%%%%%%%%%%%%%%%%%%%%%%%%%%%%%%%%%%

%%%%%%%%%%%%%%%%%%%%%%%%%%%%%%%%%%%%%%%%%%%%%%%%%%%%%%%%%%
%%%%%%%%%%%%%%%%%%%%%%%%%%%%%%%%%%%%%%%%%%%%%%%%%%%%%%%%%%
\section{Main results}\label{Bernoulli}
%%%%%%%%%%%%%%%%%%%%%%%%%%%%%%%%%%%%%%
%%%%%%%%%%%%%%%%%%%%%%%%%%%%%%%%%%%%%%

\subsection{Examples of quasi-Bernoulli independent measures}

It is not difficult to show that, in the setting of \cite{KAKI}, the
two following examples can be seen as random Gibbs measures associated
with a random H\"older potential in the dynamical system
$(\mathbb{A},\sigma)$.
%%%%%%%%%%%%%%%%%%%%%%%%%%%%%%%%%%%%%%%%%%%%%%%%%%%%%%%%%%
%%%%%%%%%%%%%%%%%%%%%%%%%%%%%%%%%%%%%%%%%%%%%%%%%%%%%%%%%%
\smallskip

\noindent
{\it Example 1. Multinomial random measures.} Let
$(W_0,\dots,W_{b-1})$ be a positive random vector such that
$\sum_{k=0}^{b-1}W_j=1$ almost surely, and let $\big
((W_0,\dots,W_{b-1})(j)\big )_{j\geq 1}$ be a sequence of independent
copies of the vector $(W_0,\dots,W_{b-1})$. Let $\ell$ denote the
unique measure on $\mathbb{A}$ such that $\ell ([w])=b^{-n}$ for $w\in
\mathcal{A}^n$.

With probability one, the sequence of measures $(\mu_j)_{j\ge 1}$
defined on $\mathbb{A}$ by
\begin{equation}
\label{muj1}
\frac{d\mu_j}{d\ell}(t)=b^j\prod_{k=1}^j W_{w_k}(k)\quad (t\in
[w_1\dots w_j])
\end{equation}
converges weakly, as $j\to\infty$, to a probability measure $\mu$
which clearly satisfies {\bf (P1)} to {\bf (P4)}. Here $\mu^{(j)}$ is
constructed like $\mu$, but with the vectors $\big
((W_0,\dots,W_{b-1})(k)\big )_{k \ge j+1}$ instead of $\big
((W_0,\dots,W_{b-1})(k)\big )_{k \ge 1}$. 
\medskip

%%%%%%%%%%%%%%%%%%%%%%%%%%%%%%%%%%%%%%%%%%%%%%%%%%%%%%%%%%
%%%%%%%%%%%%%%%%%%%%%%%%%%%%%%%%%%%%%%%%%%%%%%%%%%%%%%%%%%
\noindent
{\it Example 2. Random Riesz products.}  Let $\phi$ be a 1-periodic
  H\"older conti\-nuous function on $\mathbb{R}$ and let
  $(\theta_k)_{k\ge 0}$ be a sequence of independent random variables
  uniformly distributed in $[0,1]$. Let $\pi:\mathbb{A}\to [0,1]$ be
  the mapping $t=t_1\cdots t_k\cdots \mapsto \sum_{k\ge
  1}t_kb^{-k}$. Then consider on $\mathbb{A}$ the sequence of measures
  $(\mu_j)_{j\geq 0}$ whose density with respect to $\ell$ is given by
\begin{equation}
\label{muj2}
\frac{d\mu_j}{d\ell}(t)=\frac{\prod_{k=0}^{j-1}\exp\big (\phi(b^k\pi
(t)+\theta_k)\big )}{\int_0^1 \prod_{k=0}^{j-1}\exp\big
(\phi(b^k\pi(u)+\theta_k)\big )\, du}.
\end{equation}
Because of Theorems 3.1 and 3.2 in \cite{KAKI}, with probability one,
the sequence $\{\mu_j\}$ converges weakly to a probability measure
$\mu$. Moreover, it is shown in \cite{FAN2,BCM} that, because of the
H\"older regularity and the 1-periodicity of $\phi$, properties {\bf
(P1)} to {\bf (P3)} hold. Property {\bf (P4)} follows from the fact
that the $\theta_k$'s are chosen independent. Here $\mu^{(j)}$ is
constructed like $\mu$, but with the phases $(\theta_k )_{k\ge j+1}$
instead of $(\theta_k )_{k\ge 1}$.

%%%%%%%%%%%%%%%%%%%%%%%%%%%%%%%%%%%%%%
%%%%%%%%%%%%%%%%%%%%%%%%%%%%%%%%%%%%%%
\subsection{Identification of the function $\tau_\mu$ and auxiliary
  measures} $\ $

Let $\mu$ be quasi-Bernoulli independent random measure. We specify
the scaling function $\tau_\mu$ and the family of analysing measures
discussed in the Introduction.

$\bullet$ {\it The function $\tau_\mu$.} For every $j,k\ge 1$, let us
define the function
$$ \tau^{(k)}_j: q\in \mathbb{R} \mapsto -\frac{1}{j}\log_b\,
\sum_{w\in \mathcal{A}^j}(\mu^{(k)})_j ([w])^q,
$$ where $(\mu^{(k)})_j$ denotes the measure associated with
$\mu^{(k)}$ like $\mu_j$ is associated with $\mu$ in formulas
(\ref{muj1}) and (\ref{muj2}). When $k=0$ we simply write $\tau_j(q)$.

The same arguments as those used in \cite{FAN2} and \cite{BCM} (mainly
based on Kingman's sub-multiplicative ergodic theorem) show that, with
probability one, for all $q\in\mathbb{R}$ and for all $k\ge 0$,
$\tau^{(k)}_j(q)$ converges, as $j\to+\infty$, to a real number
$\tau_\mu(q)$ (thus independent of $k$). $\tau_\mu(q)$ coincides with
the number defined in (\ref{tau}). Moreover, $\tau_\mu(q)$ is also the
limit when $j\to+\infty$ of the sequence $\mathbb{E}\big
(\tau_j(q)\big )$. In particular the mapping $q\mapsto \tau_\mu(q)$ is
deterministic.

Due the concavity of $\tau_j$, with probability one, $\tau_j$
converges uniformly to $\tau_\mu$ on compact sets.

\smallskip 

%%%%%%%%%%%%%%%%%%%%%%%%%%%%%%%%%%%%%%
%%%%%%%%%%%%%%%%%%%%%%%%%%%%%%%%%%%%%%

$\bullet$ {\it Auxiliary measures.} The multifractal spectrum of $\mu$
is obtained thanks to the following auxiliary measures $\mu_q$. Let
$\Omega^*$ be a subset of $\Omega$ with $\mathbb{P}(\Omega^*)=1$ such
that the conclusions of Proposition ~\ref{approx} hold for all
$\omega\in\Omega^*$. For every $\omega\in \Omega^*$, for all $q\in
\mathbb{R}$ and for all $j\ge 1$, let $\mu_{q,j}$ be the probability
measure with a density with respect to the measure $\ell$ on $[v]$
(for every $v\in \mathcal{A}^j$) given by $b^{j}\mu([v])^{q}
b^{j\tau_j(q)}$.

If $\omega$ is still fixed, for every $q\in \mathbb{R}$ one can
consider a subsequence $j_n(q)$ such that the sequence
$\{\mu_{q,j_n(q)}\}_n$ converges weakly to a measure $\mu_q$ (which
depends on~$\omega$). This can also be done for the measures
$\mu^{(j)} $. For every fixed $\omega\in\Omega^*$, for all $j\ge 1$
and $q\in\mathbb{R}$, a measure $\mu^{(j)}_q$ is built as $\mu_q$.

%%%%%%%%%%%%%%%%%%%%%%%%%%%%%%%%%%%%%%

%%%%%%%%%%%%%%%%%%%%%%%%%%%%%%%%%%%%%
%%%%%%%%%%%%%%%%%%%%%%%%%%%%%%%%%%%%%
%%%%%%%%%%%%%%%%%%%%%%%%%%%%%%%%%%%%%
%%%%%%%%%%%%%%%%%%%%%%%%%%%%%%%%%%%%%
\subsection{Main results}
In the sequel, $[x]$ stands for the integer part of the real
number~$x$.  If the function $\tau_\mu$ is differentiable, $J$ stands
for the open interval $\{q\in\mathbb{R}:
\tau_\mu'(q)q-\tau_\mu(q)>0\}$.

%%%%%%%%%%%%%%%%%%%%%%%%%%%%%%%%%%%%%%
%%%%%%%%%%%%%%%%%%%%%%%%%%%%%%%%%%%%%%
\begin{theorem}\label{quasi-Bernoulli_random}
Let $\mu$ be a quasi-Bernoulli independent random measure, and assume
that $\tau_\mu$ is twice continuously differentiable. Let
$\widetilde\varepsilon =(\varepsilon_n)_{n\ge 1}$ a sequence of
positive numbers going to 0. Assume that $\forall \, (M,\alpha)>0$ the
series $\sum_{n\ge 1} b^{Mn^{3/4}\log
(n)}{b}^{-n\alpha\varepsilon_n^2}$ converges.

With probability one, $\forall \,q\in J$, the singularity sets
 $E^{\mu_q}_{\tau_\mu'(q)q-\tau_\mu(q)}(N,\widetilde \varepsilon)$ and
 $E^{\mu}_{\tau_\mu'(q)}(N,\widetilde\varepsilon)$ (defined in
 (\ref{defebm})) are both of full $\mu_q$-measure.
\end{theorem}
%%%%%%%%%%%%%%%%%%%%%%%%%%%%%%%%%%%%%%
%%%%%%%%%%%%%%%%%%%%%%%%%%%%%%%%%%%%%%
\begin{remark}\label{usa} (1) As soon as $\varepsilon_n\ge
n^{-1/8}\log(n)^{1/2+\eta}$ for some $\eta>0$, one has $b^{Mn^{ 3/4
}\log (n)}{b}^{-n\alpha\varepsilon_n^2} \leq~n^{-(1+2\eta)}$ for all
$M>0$. The conclusions of Theorem~\ref{quasi-Bernoulli_random} thus
hold in this case. 

In view of the law of the iterated logarithm (see \cite{PHI,KI}), one
could expect $\widetilde \varepsilon$ to decrease faster toward
0. This is not the case because we impose the control of neighboring
cylinders (in the sense of $\delta$) and the uniform control over the
parameter $q$.

(2) In Examples 1 and 2, $\tau_\mu$ is analytic (see \cite{BCM}
 and references therein).
\end{remark}
%%%%%%%%%%%%%%%%%%%%%%%%%%%%%%%%%%%%%%
%%%%%%%%%%%%%%%%%%%%%%%%%%%%%%%%%%%%%%
The next statement uses the definitions introduced in
Section~\ref{princ2}. The measures $\mu^{(j)}$ and $\mu^{(j)}_q$ play
respectively the role of $\mu^{(j)}_\lambda$ and $m^{(j)}_\lambda$ for
$j\ge 1$.
%%%%%%%%%%%%%%%%%%%%%%%%%%%%%%%%%%%%%%
%%%%%%%%%%%%%%%%%%%%%%%%%%%%%%%%%%%%%%
\begin{theorem}[{\bf Growth speed in singularity sets}]
\label{quasi-Bernoulli_random_ren} Under the assumptions
of Theorem~\ref{quasi-Bernoulli_random}, let us choose $\eta>0$,
$N\geq 1$ and a sequence $\widetilde\ep=(\ep_n)$ so that $\varepsilon_n \ge
n ^{-1/8}\log(n )^{1/2+\eta}$. Let us also fix $\alpha>1$.

For every compact subinterval $K$ of $J$, with probability one, for
$j$ large enough and for all $q\in K$, if $\mathcal{S}_j=\big [
\exp\big (\sqrt{\alpha \log (j)}\big )\big ]$, one has
$$\max \Big (GS\big(\mu_q^{(j)},\mu^{(j)},\tau_\mu'(q),
N,\widetilde\varepsilon \big
),GS\big(\mu_q^{(j)},\mu_q^{(j)}, \tau_\mu'(q)q-\tau_\mu(q),
N,\widetilde\varepsilon\big )\Big )\le
\mathcal{S}_j.$$
\end{theorem}
%%%%%%%%%%%%%%%%%%%%%%%%%%%%%%%%%%%%%%
%%%%%%%%%%%%%%%%%%%%%%%%%%%%%%%%%%%%%

\begin{remark}
Instead of a fixed number of neighbors, it is not difficult to treat
the case of an increasing sequence of neighbors $N_n$, simultaneously
with the speed of convergence $\varepsilon_n$. This number $N_n$ can
then go to $\infty$ under the condition that $\log
N_n=o(n\varepsilon_n^2)$.

Another improvement consists in replacing the fixed fraction $f$ in
(\ref{GS}) by a fraction $f_j$ going to 1 as $j$ goes to $\infty$. The
choice $f_j=1-b^{-s_j}$ with $s_j=o\left (\big [ \exp\big
(\sqrt{\alpha \log (j)}\big )\big ]\right )$ is convenient, as the
reader can check.
\end{remark}

Let us recall that for all integers $j\ge 1$ and $n\geq 1$
\begin{equation}
\label{defnld}
\mathcal{N}_n(\mu^{(j)},\alpha,\varepsilon_n)\!= \!\#\big\{w\in\mathcal{A}^n
  : b^{-n(\alpha+\varepsilon_n)}\le \mu^{(j)} ([w])\le
  b^{-n(\alpha-\varepsilon_n)}\big\}.
\end{equation}
%%%%%%%%%%%%%%%%%%%%%%%%%%%%%%%%%%%%%%
%%%%%%%%%%%%%%%%%%%%%%%%%%%%%%%%%%%%%%
\begin{theorem}[{\bf Speed of renewal of large deviation
spectrum}]\label{quasi-Bernoulli_random_dev} Under the assumptions of
Theorem~\ref{quasi-Bernoulli_random}, let us choose $\varepsilon_n\ge
n^{-1/8}\log(n)^{1/2+\eta}$ for some $\eta>0$. Let $K$
be a compact subinterval of $J$, and let $\beta=1+\max_{q\in
K}|q|$.

For every $\alpha>1$, with probability one, for $j$ large enough, for
all $q\in K$ and for all $n\ge \big [\exp\big (\sqrt{\alpha\log
(j)}\big )\big ]$, one has
$$
b^{n(\tau_\mu'(q)q-\tau_\mu(q)-\beta\varepsilon_n)}\le \mathcal N_n\big
(\mu^{(j)},\tau_\mu'(q),\varepsilon_n\big )\le
b^{n(\tau_\mu'(q)q-\tau_\mu(q)+\beta\varepsilon_n)}.
$$
\end{theorem}
%%%%%%%%%%%%%%%%%%%%%%%%%%%%%%%%%%%%%%
%%%%%%%%%%%%%%%%%%%%%%%%%%%%%%%%%%%%%%

The following Propositions are useful in the sequel.

%%%%%%%%%%%%%%%%%%%%%%%%%%%%%%%%%%%%%%
%%%%%%%%%%%%%%%%%%%%%%%%%%%%%%%%%%%%%%
\begin{proposition}\label{approx}
Let $K$ be a compact subset of $\mathbb{R}$, and let us fix
$\alpha>1$.

\noindent {\bf 1.}  There exists a constant $C_K$ such that 
$$ \mbox{for every }n\geq 1,\ \sup_{q\in K}\left |\mathbb{E}\big
(\tau_n(q)\big )-\tau_\mu(q)\right |\le C_K\, n^{-1}.
$$
\noindent {\bf 2.}  There exists a constant $C_K$ such that
with probability one
\begin{equation*}
\mbox{for every } n \mbox{ large enough}, \ \sup_{q\in K}\,
|\tau_n(q)-\tau_\mu(q)|\le C_K{\log(n )}{n^{-1/4}},
\end{equation*}
%%%%%%%%%%%%%%%%%%%%%%%%%%%%%%%%%%
and for $j$ large enough, for every $n\ge \big [\exp\big
(\sqrt{\alpha\log (j)}\big )\big ] $,
\begin{equation*}
\sup_{q\in K}\, |\tau^{(j)}_n(q)-\tau_\mu(q)|\le C_K{\log(n
 )}{n^{-1/4}}.
\end{equation*}
%%%%%%%%%%%%%%%%%%%%%%%%%%%%%%%%%%
\end{proposition}
%%%%%%%%%%%%%%%%%%%%%%%%%%%%%%%%%%%%%%
%%%%%%%%%%%%%%%%%%%%%%%%%%%%%%%%%%%%%%
%%%%%%%%%%%%%%%%%%%%%%%%%%%%%%%%%%%%%%

%%%%%%%%%%%%%%%%%%%%%%%%%%%%%%%%%%%%%%
%%%%%%%%%%%%%%%%%%%%%%%%%%%%%%%%%%%%%%
%%%%%%%%%%%%%%%%%%%%%%%%%%%%%%%%%%%%%%
\begin{proposition}\label{approx8}
 Assume that $\tau_\mu$ is differentiable and that $K\subset
\{q\in\mathbb{R}: \tau_\mu'(q)q-\tau_\mu (q)>0\}$. Let us denote $g_k$
the word consisting of $k$ consecutive zeros and $d_k$ the word
consisting of $k$ consecutive $b-1$.

\noindent There are three constants $(C,\eta_0,\Lambda)\in
\mathbb{R}_+^{*3}$ such that with probability one,
$$\sup_{\substack{q\in K, n\ge
0\\\gamma\in\{-1,1\},\ \eta\in(0,\eta_0]}}\left (\mu_q([d_n])
\frac{\mu([g_n])^{\gamma\eta}
}{\mu([d_n])^{\gamma\eta}}+\mu_q([g_n])
\frac{\mu([d_n])^{\gamma\eta} }{\mu([g_n])^{\gamma\eta}}\right
)b^{n \Lambda}\le C,
$$
and for $j$ large enough, for every $n\ge \big [\exp\big
(\sqrt{\alpha\log (j)}\big )\big ]$,
$$ \sup_{\substack{q\in K,\\ \gamma\in\{-1,1\},\ \eta\in(0,\eta_0]}}
\left (\mu^{(j)}_q([d_n]) \frac{\mu^{(j)}([g_n])^{\gamma\eta}
}{\mu^{(j)}([d_n])^{\gamma\eta}}+\mu^{(j)}_q([g_n])
\frac{\mu^{(j)}([d_n])^{\gamma\eta}
}{\mu^{(j)}([g_n])^{\gamma\eta}}\right)b^{n\Lambda}\le C.
$$
\end{proposition}
%%%%%%%%%%%%%%%%%%%%%%%%%%%%%%%%%%%%%%
%%%%%%%%%%%%%%%%%%%%%%%%%%%%%%%%%%%%%%
\noindent Propositions \ref{approx} and \ref{approx8} are proved in
Section \ref{secproof1}, and the theorems in Section~\ref{secproof2}.

%%%%%%%%%%%%%%%%%%%%%%%%%%%%%%%%%%%%%%
%%%%%%%%%%%%%%%%%%%%%%%%%%%%%%%%%%%%%%
\section{Proofs of Proposition \ref{approx} and \ref{approx8}}
\label{secproof1}

%%%%%%%%%%%%%%%%%%%%%%%%%%%%%%%%%%%%%%

%%%%%%%%%%%%%%%%%%%%%%%%%%%%%%%%%%%%%%%%%

%%%%%%%%%%%%%%%%%%%%%%%%%%%%%%%%%%%%%
%%%%%%%%%%%%%%%%%%%%%%%%%%%%%%%%%%%%%
\subsection{Proof of Proposition~\ref{approx}} 
\label{propapprox}

{\bf 1.} The arguments are standard. For $q\in\mathbb{R}$ and $j\ge
1$, let us define $L_j(q)=j\mathbb{E}\big (\tau_j(q)\big )$. As a
consequence of {\bf (P1)},
\begin{eqnarray*}
C^{-|q|}\sum_{v\in{\mathcal A}^j,w\in{\mathcal A}^n} (\mu_j([v])
\mu^{(j)}([w]))^q & \leq &  \sum_{v\in{\mathcal A}^{j+n}} \mu([v])^q \\
\mbox{and }\sum_{v\in{\mathcal A}^{j+n}} \mu([v])^q  & \leq  & 
C^{|q|} \sum_{v\in{\mathcal A}^j,w\in{\mathcal A}^n} (\mu_j([v])
\mu^{(j)}([w]))^q .
\end{eqnarray*}
Using then {\bf (P3)} and {\bf (P4)}, and the definition of
$\tau_j(q)$, one gets
$$ \forall \,j,n\ge 1,\ \forall \, q\in\mathbb{R}, \,\,
|L_{j+n}(q)-L_j(q)-L_n(q)|\le C_q:=|q|\log_b (C).
$$ It follows that the two sequences $L_j(q)+C_q$ and $-L_j(q)+C_q$
are sub-additive. Consequently, the sequence $(L_j(q)+C_q)/j$
converges, as $j\to + \infty$, to its infimum denoted by
$L(q)$. Similarly, the sequence $(-L_j(q)+C_q)/j$ converges to
$-L(q)$. This yields that
\begin{equation}
\label{item1}
 \forall \,j\ge 1,\,\, \forall\, q \in\mathbb{R}, \,\,\left
|{L_j(q)}{/j}-L(q)\right |\le  {C_q}{/j},
\end{equation}
which gives the desired conclusion since we have seen that
$L(q)=\tau_\mu(q)$.
 
\medskip

\noindent
{\bf 2.}  We invoke a property which does hold because of {\bf (P2)}:
there exists $M>0$ such that with probability one,
\begin{equation}\label{equicontinuite}
 \forall\, q,q' \in\mathbb{R}^2, \,\,\forall \,j\ge
1,\,\,|\tau_j(q)-\tau_j(q')|\le M |q-q'|.
\end{equation}

Fix $K$, a non-trivial compact subinterval of $\mathbb{R}$.

For $q\in K$, $j\ge 0$ and $n\ge 1$, let us define the random variables
%%%%%%%%%%%%%%
$$
L^{(j)}_n(q)=-\log_b\, \sum_{w\in \mathcal{A}^n}(\mu^{(j)})_n ([w])^q,
$$
%%%%%%%%%%%%%%
and $L_n(q)=L^{(0)}_n(q)$.  It follows from {\bf (P1)} and {\bf (P4)}
that
%%%%%%%%%%%%%%
\begin{equation}\label{subadditiveL}
\forall \ q\in K,\ \forall \, j,n\ge 1,\, \forall \, q\in K, \,
 |L_{j+n}(q)-L_j(q)-L^{(j)}_n(q)|\le |q|\log_b (C).
\end{equation}
%%%%%%%%%%%%%%
Let $C_K=\sup_{q\in K}|q|\log_b (C)$, and fix $q\in K$.  For every
integer $m\ge 1$, we write $m=[\sqrt{m}]^2+i_m$ where
$i_m\in[0,3\sqrt{m}]$. Using again {\bf (P1)} to {\bf (P4)}, one
deduces from (\ref{subadditiveL}) that for every $m\ge 1$, there exist
$[\sqrt{m}]$ independent copies $X^{(m)}_1,\dots,
X^{(m)}_{[\sqrt{m}]}$ of $L_{[\sqrt{m}]}(q)$ such that
%%%%%%%%%%%%%%
\begin{equation}\label{1}
\Big |L_m(q)-\sum_{i=1}^{[\sqrt{m}]}X^{(m)}_i(q)\Big|\le
C_K[\sqrt{m}]+|L_{i_m}(q)|\leq 4C_K\sqrt{m}.
\end{equation}
%%%%%%%%%%%%%%

We invoke the following concentration inequality (see Lemma 1.5 of
\cite{LEDOUXTALAGRAND})
%%%%%%%%%%%%%%%%%%%%%%%%%%%%%%%%%%%%%
%%%%%%%%%%%%%%%%%%%%%%%%%%%%%%%%%%%%%
\begin{lemma}\label{sum}
Let $n\ge 1$ and let $(Y_i)_{1\le i\le n}$ be a sequence of random
variable i.i.d. with a centered and bounded random variable $Y$. For
all $s>0$,
$$
\mathbb{P}\Big (\Big |\sum_{i=1}^n Y_i\Big |>\Vert Y\Vert_\infty
s\sqrt{n}\Big)\le 2 \exp \left(-{s^2}{/2}\right).
$$
\end{lemma}
%%%%%%%%%%%%%%%%%%%%%%%%%%%%%%%%%%%%%
%%%%%%%%%%%%%%%%%%%%%%%%%%%%%%%%%%%%%
Let us define the random variables
$Y^{(m)}_i(q)=X^{(m)}_i(q)-\mathbb{E}\big (X^{(m)}_i(q)\big )$. By
{\bf (P2)}, one can find a constant $M_K>0$ such that $\sup_{q\in
K}|Y^{m}_i(q)|\le M_K[\sqrt{m}]$. As a consequence, Lemma \ref{sum}
can be applied to the bounded family
$Y^{(m)}_1(q),\dots,Y^{(m)}_{[\sqrt{m}]}(q)$. Then choosing
$s=\sqrt{2}\log (m)$ yields (remember that $\|Y\|_\infty \leq
M_K[\sqrt{m}]$)
$$
\mathbb{P}\Big (\Big  |\sum_{i=1}^{[\sqrt{m}]} Y^{(m)}_i(q)\Big |>
\sqrt{2}M_K\log (m) [\sqrt{m}]^{3/2}\Big )\le \exp\big (-(\log
m)^2\big ).
$$ For every $m\ge 1$, let $q^{(m)}_1 <q^{(m)}_2<\dots
<q_k^{(m)}<\dots$ be a finite sequence of points of $K$ such that
$q^{(m)}_{k+1}-q^{(m)}_k\le m^{-1/4}$, and denote by $\mathcal {R}_m$
the set of these points. We can assume that the cardinality of
$\mathcal{R}_m$ is less than or equal to $|K|\sqrt{m+1}$. Then
\begin{eqnarray*}
&&\sum_{m\ge 1}\mathbb{P}\Big (\exists\ q\in \mathcal{R}_m,\ \Big
|\sum_{i=1}^{[\sqrt{m}]} Y^{(m)}_i(q)\Big|> \sqrt{2}M_K\log (m)
[\sqrt{m}]^{3/2}\Big )\\
&\le& \sum_{m\ge 1}|K|\sqrt{m+1}\exp\big (-(\log
m)^2\big )<\infty .
\end{eqnarray*}
This implies that for every $ q\in \mathcal{R}_m$ and for $m$ large
enough, 
\begin{equation}
\label{2}
\Big |\sum_{i=1}^{[\sqrt{m}]} Y^{(m)}_i(q)\Big| \le \sqrt{2}\log m
M_K[\sqrt{m}]^{3/2}.
\end{equation}
On the other hand, remembering the proof of item {\bf 1.} and
(\ref{item1}), one has
\begin{equation}
\label{3} \forall \,q\in \mathcal{R}_m, \ \left |\mathbb{E}\big
(X_1^{m}(q)\big )-[\sqrt{m}]\tau_\mu(q)\right|\le C_K.
\end{equation}
For every $q\in \mathcal{R}_m$, 
$\left|L_m(q)-m\tau_\mu(q)\right|$ can be upper bounded by 
\begin{eqnarray*}
\Big
|L_m(q)-\sum_{i=1}^{[\sqrt{m}]}X_i^{(m)}(q)\Big |+\Big
|\sum_{i=1}^{[\sqrt{m}]}Y^{(m)}_i(q)\Big | \\+\Big
|\sum_{i=1}^{[\sqrt{m}]}\mathbb{E}\big(X_i^{(m)}(q)\big) -[\sqrt
m]^2\tau_\mu(q)\Big |+i_{m}|\tau_\mu(q)|.
\end{eqnarray*}
With probability one, for $m$ large enough, using respectively
(\ref{1}) and (\ref{3}), this first and the third term are both bounded
by a $O \big ([\sqrt{m}]\big)$ (which does not depend on $q$). Using 
(\ref{2}) and remarking that $i_m =O \big ([\sqrt{m}]\big)$,
one gets
\begin{eqnarray*}
\left|L_m(q)-m\tau_\mu(q)\right|&\le& \sqrt{2}M_K\log(m)
[\sqrt{m}]^{3/2}+O \big ([\sqrt{m}]\big),
\end{eqnarray*}
where $O\big ([\sqrt{m}]\big)$ is uniform over $q\in
\mathcal{R}_m$. This yields $ |\tau_m(q)-\tau_\mu(q)|=O\Big
({\log (m)}{m^{-1/4}}\Big ) $ uniformly for $q\in
\mathcal{R}_m$ when $m$ is large enough.
The  conclusion follows from (\ref{equicontinuite}) and from the
construction of the sets $\mathcal{R}_m$.

\smallskip

Let us show the second inequality of item {\bf 2.}  For every $j\ge
0$, $m\ge 1$ and $q\in K$, let us consider a sequence
$Y^{(m),j}_i(q)$, $1\le i\le [\sqrt{m}]$, associated with
$\mu=\mu^{(j)}$ like $Y^{(m)}_i(q)$, $1\le i\le [\sqrt{m}]$, is
associated with $\mu=\mu^{(0)}$. Let $\mathcal{R}_m$ be defined as
above, and let us consider the events
$$
A(j,m)=\Big \{\exists\ q\in \mathcal{R}_m,\ \Big
|\sum_{i=1}^{[\sqrt{m}]} Y^{(m),j}_i(q)\Big|> \sqrt{2}M_K\log (m)
[\sqrt{m}]^{3/2}\Big \}.
$$ One verifies that $\sum_{j\ge 0}\sum_{m\ge \,  [\, \exp \sqrt{\alpha
\log (j)}  \,]\,}\mathbb{P}\big (A(j,m)\big)<\infty$. We then deduce
from the Borel-Cantelli Lemma that with probability one, for $j$ large
enough, if $m\ge \big [\exp\big (\sqrt{\alpha \log (j)}\big )\big ]$
then $A(j,m)^c$ holds. One concludes by using the same estimates as
above.

%%%%%%%%%%%%%%%%%%%%%%%%%%%%%%%%%%%%%
%%%%%%%%%%%%%%%%%%%%%%%%%%%%%%%%%%%%%
\subsection{Proof of Proposition~\ref{approx8}} 
\label{propapprox8}

If $t_j\in\{g_j,d_j\}$, the same kind of arguments as in the
proof of Proposition~\ref{approx} show that, with probability one,
$\Lambda_t(1)=\lim_{j\to\infty}\frac{1}{j}\log_b \big
(\mu([t_j])\big )$ exists, and this number is deterministic. Hence,
using (\ref{quasibmuq}), with probability one, for every $q\in
\mathbb{R}$, the limit
$\Lambda_t(q)=\lim_{j\to\infty}\frac{1}{j}\log_b \big
(\mu_q([t_j])\big )$ exists and is equal to $q\Lambda_t(1)+\tau_\mu
(q)$. Since $\mu_q$ is a finite measure, $\Lambda_t(q)\le 0$.

Moreover, there exists $C_K>0$ such that for $\mathbb{P}$-almost every
$\omega\in\Omega^*$, for $j$ large enough, for all $q\in K\cup\{1\}$,
\begin{eqnarray*}
&&\Big |\frac{1}{j}\log_b \big (\mu_q([t_j])\big
)-\Lambda_t(q)\Big|\le C_K{\log (j)}{j^{-1/4}}\\
\forall \,k\ge \big [\exp\big
(\sqrt{\alpha \log (j)}\big )\big ], &&\Big |\frac{1}{k}\log_b \big
(\mu^{j}_q([t_k])\big )-\Lambda_t(q)\Big |\le C_K{\log
(k)}{k^{-1/4}}.
\end{eqnarray*}
So, for $j$ large enough, $\gamma\in
\{-1,1\}$ and $\eta>0$, one has
%%%%%%%%%%%%%%%%%%%%%%%%%%%%%%%%%%%%%%%
%%%%%%%%%%%%%%%%%%%%%%%%%%%%%%%%%%%%%%%
\begin{eqnarray*}
&&\mu_q([d_j]) \frac{\mu([g_j])^{\gamma\eta}
}{\mu([d_j])^{\gamma\eta}}+\mu_q([g_j]) \frac{\mu([d_j])^{\gamma\eta}
}{\mu([g_j])^{\gamma\eta}} \le f(j)
\end{eqnarray*}
and $ \forall \,k\ge \big [\exp\big (\sqrt{\alpha \log (j)}\big )\big
]$
$$
\mu^{(j)}_q([d_k]) \frac{\mu^{(j)}([g_k])^{\gamma\eta}
}{\mu^{(j)}([d_k])^{\gamma\eta}}+\mu^{(j)}_q([g_k])
\frac{\mu^{(j)}([d_k])^{\gamma\eta}
}{\mu^{(j)}([g_k])^{\gamma\eta}} \le f(k),
$$
%%%%%%%%%%%%%%%%%%%%%%%%%%%%%%%%%%%%%%%
%%%%%%%%%%%%%%%%%%%%%%%%%%%%%%%%%%%%%%%
where
$$
f(j)=b^{(2\eta +1)C_Kj^{3/4}\log (j)}b^{j\eta \big (|\Lambda_d
(1)|+|\Lambda_g (1)|\big )}\big
(b^{j\Lambda_d(q)}+b^{j\Lambda_g(q)}\big ).
$$
%%%%%%%%%%%%%%%%%%%%%%%%%%%%%%%%%%%%%%%%%%%
%%%%%%%%%%%%%%%%%%%%%%%%%%%%%%%%%%%%%%%%%%%
Let us show that $\Lambda_t(q)<0$ for $t\in\{g,d\}$ and $q\in
K$. Suppose $\Lambda_t(q_0)=0$ for some $q_0\in K$. Remember that
$q\Lambda_t(1)+\tau_\mu (q)=\Lambda_t(q)\le 0$ for all
$q\in\mathbb{R}$. Using the concavity of $\tau_\mu$, the equality
$\Lambda_t(q_0)=0$ implies that $\tau'_\mu (q_0)=-\Lambda_t(1)$ and
then that $\tau_\mu'(q_0)q_0-\tau_\mu(q_0)=0$, in contradiction with
our assumption $K\subset J$.

Finally, since $\Lambda_t(q)=q\Lambda_t(1)+\tau_\mu (q)$, the mapping
$q\mapsto \Lambda_t(q)$ is continuous, and the conclusion follows from
properties {\bf (P1)} and {\bf (P2)}, the compactness of $K$, and the
form of $f(j)$.
 
%%%%%%%%%%%%%%%%%%%%%%%%%%%%%%%%%%%%%
%%%%%%%%%%%%%%%%%%%%%%%%%%%%%%%%%%%%%
%%%%%%%%%%%%%%%%%%%%%%%%%%%%%%%%%%%%%
%%%%%%%%%%%%%%%%%%%%%%%%%%%%%%%%%%%%%

%%%%%%%%%%%%%%%%%%%%%%%%%%%%%%%%%%%%%%
%%%%%%%%%%%%%%%%%%%%%%%%%%%%%%%%%%%%%%
\section{Proofs of Theorems \ref{quasi-Bernoulli_random},  \ref{quasi-Bernoulli_random_ren} and \ref{quasi-Bernoulli_random_dev} } 
\label{secproof2}

Mimicking the approach in \cite{BCM} and using Proposition
\ref{approx} and \ref{approx8} shows that if $K$ is a compact subset
of $\mathbb{R}$, there exists a constant $M_K$ such that, with
probability one, for $n$ large enough, $\forall \,v\in
\mathcal{A}^{n}$, $\forall \,q\in K$,
\begin{equation}\label{quasibmuq}
M_K^{-1}b^{-(M_K) n^{3/4}\log (n)}\le
\frac{\mu_q([v])}{\mu([v])^q b^{n\tau_\mu(q)}}\le M_Kb^{(M_K)n^{
3/4}\log (n)}.
\end{equation}
and if $\alpha >1$ is fixed, for the same constant $M_K$, with
probability one, for $j$ large enough, for $n\ge\big [\exp\big
(\sqrt{\alpha\log(j)}\big ) \big ]$, $\forall \,v\in \mathcal{A}^{n}$,
$\forall \,q\in K$,
\begin{equation}\label{quasibmuq2}
 M_K^{-1}b^{-(M_K) n^{3/4}\log (n)}\le
\frac{\mu^{(j)}_q([v])}{\mu^{(j)}([v])^q b^{n\tau_\mu(q)}}\le
M_Kb^{(M_K)n^{ 3/4}\log (n)}.
\end{equation}

Before starting the proofs, let us make a last useful remark.
%%%%%%%%%%%%%%%%%%%%%%%%%%%%%%%%%%%%%%%%%
\begin{remark}\label{reduction}
If $v$ and~$w$ are words of length~$n$, and if $\bar{v}$ and $\bar{w}$
stand for their prefixes of length~$n-1$, then
$\delta(\bar{v},\bar{w})> k$ implies $\delta(v,w) >bk$. It implies
that, given two integers $n \ge m>0$ and two words $v$ and $w$ in
${\mathcal A}^n$ such that $b^{m-1}<\delta(v,w)\le b^m$, there are two
prefixes $\bar{v}$ and $\bar{w}$ of respectively $v$ and $w$ of common
length\ $n-m$\ such that $\delta(\bar{v},\bar{w})\le 1$; moreover, for
these words $\bar{v}$ and $\bar{w}$, there are at most $b^{2m}$ pairs
$(v,w)$ of words in $\mathcal{A}^n$ such that $\bar{v}$ and $\bar{w}$
are respectively the prefixes of $v$ and $w$.
\end{remark} 
%%%%%%%%%%%%%%%%%%%%%%%%%%%%%%%%%%%%%%
%%%%%%%%%%%%%%%%%%%%%%%%%%%%%%%%%%%%%%
\subsection{Proof of Theorem \ref{quasi-Bernoulli_random}}
Fix $K$ a compact subinterval of $J$ and $( \eta_n)_{n\ge 1}$ a
bounded positive sequence to be precised later. For
$\omega\in\Omega^*$ and $q\in K$, let us introduce the two quantities
(recall (\ref{defsneej}))
\begin{equation}
\label{deffj}
F_n(q)=S_n^{N,\varepsilon_n,\eta_n}\big (\mu_q,\mu,\tau_\mu'(q)\big )
\mbox{ and } G_n(q)=S_n^{N,\varepsilon_n,\eta_n}\big
(\mu_q,\mu_q,\tau_\mu'(q)q-\tau_\mu(q)\big ).
\end{equation}

Due to Proposition~\ref{principe1}, we seek for a uniform control of
$F_n$ and $G_n$ on $K$. We only consider $F_n$, since the study of
$G_n$ is similar. 

$\bullet$ {\bf An upper bound for $ F_n(q)$:} Consider $v,w\in
\mathcal{A}^{n}$ such that $\delta(v,w)=k\le N$, as well as two
prefixes $\bar{v}$ and $\bar{w}$ of respectively $v$ and $w$ of common
length\ $n-[\log_b(k)]$ such that $\delta(\bar{v},\bar{w})\le 1$. Let
$q\in K$. If $n$ is large enough, (\ref{quasibmuq}) holds for both $v$
and $\bar{v}$. Then, using the construction of $\mu_q$, item {\bf 2.}
of Proposition~\ref{approx}, {\bf (P1)}, {\bf (P2)} and
(\ref{quasibmuq}), one gets for $n$ large enough
%%%%%%%%%%%%%%
\begin{equation}\label{hihi}
\mu_q([v])\le \widetilde Cb^{\widetilde Cn^{ 3/4 }\log
(n)}\mu_q([\bar{v}])\quad\mbox{and}\quad \mu([w])^{\gamma \eta_n}\le
\widetilde C\mu([\bar{w}])^{\gamma \eta_n},
\end{equation}
%%%%%%%%%%%%%
where $\widetilde C$ depends on $C$, $K$, $\Vert \widetilde
\eta\Vert_\infty$ and $\Vert \widetilde\varepsilon\Vert_\infty$. Thus,
by Remark \ref{reduction}, for $n$ large enough, $0\le k\le N$ and
$\gamma\in\{-1,1\}$,
\begin{eqnarray*}
\! \! &\! \! &\!  b^{n(\tau_\mu'(q)-\gamma\varepsilon_n)\gamma\eta_n}
\hspace{-4mm} \sum_{v,w\in \mathcal{A}^{n}, \,\delta(v,w)=k}
\hspace{-4mm} \mu_q([v])\mu([w])^{\gamma \eta_n}\\
\! \! &\leq\!  & \! \widetilde C b^{\widetilde C n^{ 3/4 }\log (n)}
b^{(n-[\log_b k])(\tau_\mu'(q) -\gamma\varepsilon_n)\gamma\eta_n}
\hspace{-8mm} \sum_{v,w\in \mathcal{A}^{n-[\log_b k]}, \,
\delta(v,w)\le 1} \hspace{-8mm} \mu_q([v])\mu([w])^{\gamma\eta_n},
\end{eqnarray*}
for some other constant $\widetilde C$ depending on $C$, $K$, $N$,
$\Vert \widetilde\varepsilon\Vert_\infty$ and $\Vert \widetilde
\eta\Vert_\infty$ .

Let us remark that for every integer $l \in \{0,..,\log_b (N)\}$, there
are less than $b^{l+1}$ integers $k\in [0,N]$ such that $[\log_b
k]=l$.  One thus deduces from the definition (\ref{deffj}) (and
(\ref{defsneej})) of $F_n(q)$ and from the above estimate that
\begin{equation}
\label{qwer}
F_n(q)\le\widetilde C \, b^{\widetilde C n^{ 3/4 }\log
(n)}\sum_{l=0}^{[\log_b (N)]+1}b^{l+1}
\big(T_1(q,n,l)+T_2(q,n,l)\big),
\end{equation}
where
\begin{eqnarray}
\label{deft}
 T_1(q,n,l)\!\! &=&\!\!\!\!\sum_{\gamma\in\{-1,1\}}b^{(n-l)(\tau_\mu'(q)
-\gamma\varepsilon_n)\gamma\eta_n}\sum_{w\in \mathcal{A}^{n-l}}
\mu_q([w])\mu([w])^{\gamma\eta_n}
\\ \label{deft2}
T_2(q,n,l)\!\!& =& \!\!\!\!\sum_{\gamma\in\{-1,1\}}b^{(n-l)(\tau_\mu'(q)
  -\gamma\varepsilon_n)\gamma\eta_n}\hspace{-4mm} \sum_{v,w\in
  \mathcal{A}^{n-l},\,\delta(v,w)= 1}\hspace{-5mm}\!\!
\mu_q([v])\mu([w])^{\gamma\eta_n}.\!\!\!\!\!\!\!
\!\!\!\!
\end{eqnarray}
Let us first upper bound $T_1(q,n,l)$. (\ref{quasibmuq}) yields for
some constant $M'_K$ that 
$$\sum_{w\in \mathcal{A}^{m}} \mu_q([w])\mu([w])^{\gamma\eta_n} \leq
b^{m\tau_\mu(q)} M'_K b^{M'_Km^{ 3/4 }\log m}\sum_{w\in
\mathcal{A}^{m}} \mu([w])^{q+\gamma\eta_n}, $$
where $m=n-l$. Using item {\bf 2.} of Proposition~\ref{approx}, for
some constant $C_K$
$$ \sum_{w\in \mathcal{A}^{n-l}} \mu([w])^{q+\gamma\eta_n} \leq
b^{-(n-l)\tau_\mu(q+\gamma\eta_n) +(n-l) C_K\log (n-l)
(n-l)^{-1/4}}.
$$ Since $\tau_\mu$ is twice continuously differentiable, one has
$\tau_\mu(q+\gamma\eta_n)-\tau_\mu(q)-\gamma\eta_n\tau_\mu'(q) =
\eta_n O(\eta_n)$ independently of $q\in K$ (if
$\Vert\widetilde\eta\Vert$ is small enough), and
\begin{equation}
\label{T1}
T_1(q,n,l)\le 2M'_Kb^{(M'_K+C_K)(n-l)^{ 3/4 } \log
(n-l)}b^{-(n-l) \eta_n\left (\varepsilon_n+O(\eta_n)\right )}.
\end{equation}
In order to estimate $T_2(q,n,l)$, we use the words $g_k$ and
$d_k$ defined in Proposition \ref{approx8}.  For every $m\geq 1$,
a representation of the set of pairs $(v,w)$ in ${\mathcal A}^m$ such
that $\imath(w)=\imath(v)+1$ is the following:
\begin{equation}
\label{decoupage}
\bigcup_{k=0}^{m-1} \ \bigcup_{u\in {\mathcal A}^{m-1-k}} \
\bigcup_{r\in \{0,\dots,b-2\}} \big \{(u.r.d_k,
u.(r+1).g_k)\big \}.
\end{equation}
Let $m=n-l$. The sum $\mathcal{T}_{n,\gamma}(q)=\sum_{v,w\in
\mathcal{A}^{m}, \ \delta(v,w)= 1} \mu_q([v])\mu([w])^{\gamma\eta_n}$
equals
\begin{eqnarray*}
\sum_{k=0}^{m-1} \sum_{u\in {\mathcal A}^{m-1-k}}  \sum_{r=0}^{b-2}\
\sum_{(e,f)\in\{(d,g),(g,d)\}} \mu_q([u.r.e_k])
\mu([u.(r+1).f_k])^{\gamma\eta_n}.
%+\mu([u.r.d_k])^{\gamma\eta_n}\mu_q([u.(r+1).g_k]).
\end{eqnarray*}
Let us introduce $
\Theta(q,k,n,\gamma)=\mu([d_k])^q\mu([g_k])^{\gamma\eta_n}
+\mu([g_k])^q\mu([d_k])^{\gamma\eta_n}.  $ Using (\ref{quasibmuq}) and
property {\bf (P1)} of $\mu$, one obtains another constant $\widetilde
C$ such~that
$$
\mathcal{T}_{n,\gamma}(q)\le \widetilde C (b-2)b^{C_Km^{ 3/4 }\log
(m)}b^{m\tau_\mu(q)} \sum_{k=0}^{m-1}\Theta(q,k,n,\gamma) \!\!
\sum_{u\in {\mathcal A}^{m-1-k}}\!\!  \mu([u])^{q+\gamma\eta_n} .
$$ Then, item {\bf 2.} of Proposition~\ref{approx} yields (with
another $\widetilde C$)
\begin{eqnarray*}
&&\mathcal{T}_{n,\gamma}(q) \\ & \le &\widetilde C \, b^{C_Km^{ 3/4
}\log (m)}b^{m\tau_\mu(q)}\\ && \quad\quad\times \sum_{k=0}^{m-1}
\Theta(q,k,n,\gamma) b^{2C_K(m-1-k)^{ 3/4 }\log
(m-k-1)-(m-k-1)\tau_\mu(q+\gamma\eta_n)}\\
&\le &\widetilde C b^{2C_Km^{
3/4 }\log (m)} b^{m (\tau_\mu(q)-\tau_\mu(q+\gamma\eta_n) )}
\sum_{k=0}^{m-1} \Theta(q,k,n,\gamma)
b^{(k+1)\tau_\mu(q+\gamma\eta_n)}\\
&\le & \widetilde C b^{m (-\tau_\mu'(q)\gamma\eta_n +O(\eta_n^2)
)}b^{2C_Km^{ 3/4 }\log (m)} \sum_{k=0}^{m-1}\Theta(q,k,n,\gamma)
b^{k\tau_\mu(q+O(\eta_n))}.
\end{eqnarray*}
By Proposition~\ref{approx8} and (\ref{quasibmuq}), the sum
$\sum_{k=0}^{m-1}\Theta(q,k,n,\gamma) b^{k\tau_\mu(q+O(\eta_n))}$ is
uniformly bounded over $q\in K$ and $m\ge 0$ when $\|\widetilde
\eta\|_\infty$ is small enough. Hence, replacing $m$ by $n-l$,
\begin{equation}
\label{T2}
T_2(q,n,l)\le \widetilde Cb^{2C_K (n-l)^{ 3/4 }\log (n-l)}\,
b^{-(n-l)\eta_n  (\varepsilon_n+O(\eta_n)  )}.
\end{equation}
Finally, combining (\ref{qwer}), (\ref{T1}) and (\ref{T2}) yields
\begin{eqnarray*}
F_n(q) \! \! \! \! &\le &\! \!\!\widetilde
C b^{\widetilde C n^{ 3/4 }\log (n)}\!
\!\!\sum_{l=0}^{[\log_b (N)]+1} \! \!b^{l+1}   b^{O  ((n-l)^{ 3/4 }\log (n-l) 
)}b^{-(n-l)\eta_n  (\varepsilon_n+O(\eta_n)  )}\\
&=&\!\! O \big  (b^{Mn^{ 3/4 }\log (n)}{b}^{-n\eta_n 
(\varepsilon_n+O(\eta_n)  )}  \big )
\end{eqnarray*}
for some $M>0$ independently of $q\in K$. By our assumption on
$\varepsilon_n$, the choice $\eta_n=\alpha \varepsilon_n$ with
$\alpha$ small enough so that $\eta_n (\varepsilon_n+O(\eta_n)
)\ge\alpha\varepsilon_n^2/2$ makes the series $\sum_{n\ge 1} F_n(q)$
converge for every $q\in K$. The conclusion concerning the sets
$E^\mu_{\tau_\mu'(q)}(N,\widetilde\varepsilon)$ then follows from
Proposition~\ref{principe1}.

\smallskip

%An alternative proof of the result concerning the sets
%$E^\mu_{\tau_\mu'(q)q-\tau_\mu(q)}( N,\widetilde\varepsilon)$ can be
%obtained by using (\ref{quasibmuq}) and the results for the sets
%$E^\mu_{\tau_\mu'(q)}(N,\widetilde\varepsilon)$.

%%%%%%%%%%%%%%%%%%%%%%%%%%%%%%%%%%%%%
%%%%%%%%%%%%%%%%%%%%%%%%%%%%%%%%%%%%%
\subsection{Proof of Theorem~\ref{quasi-Bernoulli_random_ren}} Fix
$\bar{\alpha}\in (1,\alpha)$. We use twice (\ref{quasibmuq2}), with
$\bar{\alpha}$ and $\alpha$, in order to get a control like
(\ref{hihi}). For $j$ large enough, if $n\ge \big[\exp\big
(\sqrt{\alpha\log(j)}\big )\big ]$ and if $v\in \mathcal{A}^n$ and
$\bar{v}$ is a prefix of $v$ such that $|\bar v|\ge n-\log_b(n)$, then
$|\bar v|\ge \big[\exp\big (\sqrt{\bar{\alpha} \log(j)}\big )\big ]$
and (\ref{hihi}) holds for $v$ and $\bar{v}$. Then, from the
computations performed in the proof of
Theorem~\ref{quasi-Bernoulli_random} and from
Proposition~\ref{approx8}, one deduces that for every compact
subinterval $K$ of $J$, there exist $C,M,\beta>0$ and $\widetilde\eta
=(\eta_n)_{n\ge 1}\in \mathbb{R}_+^{\mathbb{N}^*}$ such that with
probability one, for $j$ large enough, if $n\ge \big[\exp\big
(\sqrt{\alpha\log(j)}\big )\big ]$, for all $q\in K$

\begin{eqnarray*}
&& \max \Big (S^{N,\varepsilon_n,\eta_n}\big
(\mu_q^{(j)},\mu^{(j)},\tau_\mu'(q)\big ),
S^{N,\varepsilon_n,\eta_n}\big
(\mu_q^{(j)},\mu_q^{(j)},q\tau_\mu'(q)-\tau_\mu(q)\big )\Big ) \\ &&
\le C b^{Mn^{ 3/4 }\log (n)-\beta n\varepsilon^2_n}.
\end{eqnarray*}
In order to apply Proposition~\ref{ren}, let us define

\smallskip
\noindent
$\bullet$ $\Lambda=K$, $\lambda =q$ and
$\{(m^{(j)}_\lambda,\mu^{(j)}_\lambda)\}_{j\ge 0,\lambda \in
K}=\{(\mu_q^{(j)},\mu^{(j)})\}_{j\ge 0,q \in K}$,

\smallskip
\noindent
$\bullet$
$\{\beta_\lambda\}_{\lambda\in \Lambda}=\{\tau_\mu'(q)\}_{q\in K}$,
 
\smallskip
\noindent
$\bullet$ for every $j\ge 1$ and for $n\ge 1$, $V_n^{(j)}=C
b^{Mn^{3/4}\log (n)}b^{-\beta n\varepsilon_n^2}$,

\smallskip
\noindent
$\bullet$ for every $j\ge 1$,
$\psi_j(\widetilde\eta)=\mathcal{S}_j=\big [\exp\big
(\sqrt{\alpha\log(j)}\big )\big ]$.

\smallskip

With these parameters the conditions of Proposition~\ref{ren} are
fulfilled. As a consequence, with probability one, for $j$ large
enough, for all $q\in K$,
$GS(\mu_q^{(j)},\mu^{(j)},\tau_\mu'(q),N,\widetilde \varepsilon)\le
\big [\exp\big (\sqrt{\alpha\log(j)}\big )\big ]$.

Let us then consider the families $\{(\mu_q^{(j)},\mu_q^{(j)})\}_{j\ge
0,q \in K}$ and $\{\tau_\mu'(q)q-\tau_\mu(q)\}_{q\in K}$ instead of
the family $\{(\mu_q^{(j)},\mu^{(j)})\}_{j\ge 0,q \in K}$ and
$\{\tau_\mu'(q)\}_{q\in K}$ respectively, and keep the same
definitions for the other variables involved in
Proposition~\ref{ren}. Then the same control as above holds for the
growth speed
$GS\big(\mu_q^{(j)},\mu_q^{(j)},q\tau_\mu'(q)-\tau_\mu(q),N,\widetilde
\varepsilon\big)$. Notice that here the vector $V^{(j)}$ is
deterministic.

%%%%%%%%%%%%%%%%%%%%%%%%%%%%%%%%%%%%%
%%%%%%%%%%%%%%%%%%%%%%%%%%%%%%%%%%%%%
\subsection{Proof of  Theorem~\ref{quasi-Bernoulli_random_dev}}
Fix $\alpha>1$ and $K$ a compact subinterval of $J$ and $
N=0$. A standard Markov inequality (as in Proposition \ref{principe1})
shows that for $j\ge 0$, $n\ge 1$ and $q\in K$, one has $\mathcal
N_n\left(\mu^{(j)},\tau_\mu'(q), \varepsilon_n\right) \leq
b^{-n\tau^{(j)}_n(q)} b^{nq(\tau_\mu'(q)+sgn(q)\varepsilon_n)}$, where
$sgn(q)$ stands for the sign of $q$. Then, by
Proposition~\ref{approx8}, with probability one, one has for $j$ large
enough, for $n\ge \exp\big (\sqrt {\alpha \log (j)}\big )$ and for
$q\in K$
$$\mathcal N_n\big(\mu^{(j)},\tau_\mu'(q),
\varepsilon_n\big) \leq
b^{-n\tau^{(j)}_n(q)} b^{nq(\tau_\mu'(q)+sgn(q)\varepsilon_n)} \le
b^{n(\tau_\mu'(q)q-\tau_\mu(q)+\varepsilon'_n)},
$$where $\varepsilon'_n=\sup_{q\in K}M_Kn^{-1/4}\log
(n)+|q|\varepsilon_n$. One remarks that $\varepsilon'_n \le
\varepsilon_n (1+\sup_{q\in K}|q|) $ for $n$ large enough.  On the
other hand, let 
$$E= \Big ( E^{\mu^{(j)}}_{\tau_\mu'(q),\big [\exp\sqrt {\alpha \log (j)}\big
]}(0,\widetilde\varepsilon) \Big ) \bigcap \Big ( 
E^{\mu_q^{(j)}}_{\tau_\mu'(q)q-\tau_\mu(q), \big[\exp \sqrt {\alpha
\log (j) }\big ]}(0,\widetilde\varepsilon) \Big ).$$ Using
Theorem~\ref{quasi-Bernoulli_random_ren}, with probability one, for
$j$ large enough, for all $q\in N$, $\mu_q^{(j)} (E )\ge
{\Vert\mu_q^{(j)}\Vert}/2=1/2$. But, looking back at the definition of
$E$, one remarks that $\mu_q^{(j)} (E )\leq \mathcal N_n\big
(\mu^{(j)},\tau_\mu'(q), \varepsilon_n
\big)b^{-n(\tau_\mu'(q)q-\tau_\mu(q)-\varepsilon_n)}$ for $n\ge
\exp\big (\sqrt {\alpha \log (j)}\big )$, that is $
b^{n(\tau_\mu'(q)q-\tau_\mu(q)-\varepsilon_n)} /2 \le \mathcal
N_n\left(\mu^{(j)},\tau_\mu'(q), \varepsilon_n\right) $.

%%%%%%%%%%%%%%%%%%%%%%%%%%%%%%%%%%%%%%%%%%%%%%
%%%%%%%%%%%%%%%%%%%%%%%%%%%%%%%%%%%%%%%%%%%%%
%%%%%%%%%%%%%%%%%%%%%%%%%%%%%%%%%%%%%%%%%%%%%%
%%

%%%%%%%%%%%%%%%%%%%%%%%%%%%%%%%%%%%%%%%%%%%%%%%%
%%%%%%%%%%%%%%%%%%%%%%%%%%%%%%%%%%%%%%%%%%%%%%%%

%%%%%%%%%%%%%%%%%%%%%%%%%%%%%%%%%%%%%%%%%%%%%%%%
%%%%%%%%%%%%%%%%%%%%%%%%%%%%%%%%%%%%%%%%%%%%%%%%
%%%%%%%%%%%%%%%%%%%%%%%%%%%%%%%%%%%%%%%%%%%%%%%%

\end{document}